\documentclass[11pt]{smfart}

\usepackage{amsmath,a4wide}
\usepackage{amssymb}
\usepackage{xspace}
\usepackage{euscript}
\usepackage{mathrsfs} 
\usepackage[all]{xy}
\usepackage{hyperref}
\usepackage[francais]{babel}
\usepackage[T1]{fontenc}
\usepackage{textcomp}
\usepackage{smfthm}

\newtheorem{lemme}[subsection]{Lemme}		

\newtheorem{proposition}[subsection]{Proposition}	

\newtheorem*{thprinc}{Théorème principal}

\newtheorem{theoreme}[subsection]{Th\'eor\`eme}

\newtheorem{corollaire}[subsection]{Corollaire}

\theoremstyle{definition}

\newtheorem*{remarqueintro}{Remarque}
\newtheorem{remarque}[subsection]{Remarque}

\newtheorem{exemple}[subsection]{Exemple}

\newenvironment{demo}{\noindent{\textit{D\'emonstration. --- }}}{~\qedsymbol \vspace{4mm}}

\renewcommand{\theenumi}{\textbf{\roman{enumi}}} 

\renewcommand\theequation{\thesection.\textbf{\alph{equation}}} 

\makeatletter
\@addtoreset{equation}{subsection}
\makeatother

\renewcommand\thesubsubsection{\textbf{\thesubsection.\arabic{subsubsection}}}
\renewcommand\thesubsection{\textbf{\thesection.\arabic{subsection}}}

\newcommand{\Z}{\mathbb Z}
\newcommand{\Q}{\mathbb Q}
\newcommand{\R}{\mathbb R}

\newcommand{\F}{\mathbb F}

\newcommand{\Gm}{{\mathbb G}_m}

\newcommand{\Hom}{\underline{\mathrm{Hom}}}

\newcommand{\Spec}{\mathrm{Spec}}

\newcommand{\iso}{\buildrel{\sim}\over\rightarrow}

\renewcommand{\L}{\mathrm{L}}		

\def\commutatif{\ar@{}[rd]|{\circlearrowleft}}
\def\cartesien{\ar@{}[rd]|{\square}}

\renewcommand{\lim}{{\mathrm{lim}}} 

\def\sga#1#2#3{[{\bf $\mathbf{SGA\,{#1}}$}~{\sc #2}~#3]} 	
\def\egazéro#1#2{[{\bf \'EGA}~$0_{\textsc{#1}}$~#2]}		

\numberwithin{equation}{subsubsection}
\renewcommand{\theequation}{\thesubsection.\Alph{equation}}

\newcommand{\isolong}{\buildrel{\sim}\over\longrightarrow}
\newcommand{\pib}{\overline{\pi}}

\newcommand{\gsp}{\mathrm{GSp}_{2g}}
\newcommand{\ag}{\mathcal{A}_{g}}
\newcommand{\agn}{\mathcal{A}_{g,n}}
\newcommand{\agd}{\mathcal{A}_{g,\mathscr{D}}}
\newcommand{\agnd}{\mathcal{A}_{g,n,\mathscr{D}}}

\newcommand{\agdt}{\tilde{\mathcal{A}}_{g,\mathscr{D}}}
\newcommand{\agndt}{\tilde{\mathcal{A}}_{g,n,\mathscr{D}}}
\newcommand{\agb}{\overline{\mathcal{A}}_{g}}
\newcommand{\agnb}{\overline{\mathcal{A}}_{g,n}}
\newcommand{\agdb}{\overline{\mathcal{A}}_{g,\mathscr{D}}}
\newcommand{\agndb}{\overline{\mathcal{A}}_{g,n,\mathscr{D}}}

\newcommand{\mv}{\mathcal{M}_{V'}}
\newcommand{\mvd}{\mathcal{M}_{V',\mathscr{D}}}
\newcommand{\mvb}{\overline{\mathcal{M}}_{V'}}
\newcommand{\mvdb}{\overline{\mathcal{M}}_{V',\mathscr{D}}}
\newcommand{\bv}{\mathcal{B}_{V'}}
\newcommand{\bvd}{\mathcal{B}_{V',\mathscr{D}}}
\newcommand{\av}{\mathcal{A}_{V'}}
\newcommand{\avd}{\mathcal{A}_{V',\mathscr{D}}}

\title{Sous-groupes canoniques partiels}
\author{Vincent Pilloni et Beno\^it Stroh}
\date{20 mars 2009}
\email{pilloni@math.columbia.edu}
\email{stroh@math.univ-paris13.fr}
\address{Université Paris 13, LAGA,
99 avenue J.B. Clément,
93430 Villetaneuse
France}

\begin{document}
\maketitle

\begin{abstract} La réduction des variétés de Siegel modulo un nombre premier~$p$ est stratifiée par le rang multiplicatif du groupe~$p$-divisible de la variété abélienne universelle. Pour~$r\geq 0$, le sous-groupe multiplicatif maximal de la restriction du groupe de $p$-torsion de la variété abélienne universelle à la $r$-ième strate se relève canoniquement sur le tube de cette strate et définit un \og sous-groupe canonique partiel de rang~$r$\fg. Nous montrons qu'il existe un voisinage strict du tube sur lequel ce sous-groupe s'étend de manière finie et plate. Sur la strate ordinaire et au voisinage de celle-ci, on retrouve le sous-groupe canonique usuel étudié par Abbès et Mokrane, et Andreatta et Gasbarri.
\end{abstract}

\begin{altabstract} The reduction of Siegel varieties modulo a prime number~$p$ is stratified by the multiplicative rank of the $p$-divisible group of the universal abelian variety. For~$r\geq 0$, the maximal multiplicative subgroup of the restriction of the $p$-torsion group of the universal abelian variety to the $r$-th stratum lifts canonically to the tube of this stratum and defines a \og partial canonical subgroup of rank~$r$ \fg. We show that there exists a strict neighbourhood of the tube on which this subgroup extends in a finite flat way. On the ordinary stratum and its neighbourhood, we recover the usuel canonical subgroup studied by Abbes and Mokrane, and Andreatta and Gasbarri.
 \end{altabstract}

Soient~$1\leq r \leq g$ deux entiers, $p$ un nombre premier et $n\geq 3$ un entier non divisible par~$p$. Notons~$\agn$ le schéma sur $\Spec(\Z[1/n])$ qui paramètre les vari\'et\'es ab\'eliennes principalement polaris\'ees de genre~$g$ munies d'une structure de niveau principale en~$n$, et notons~$G$ le sch\'ema ab\'elien universel sur~$\agn$.
Soit~$\agnb$ la compactification toro\"idale de~$\agn$ associ\'ee au choix d'une d\'ecomposition poly\'edrale admissible polarisée~\cite[ch. IV]{Deg@FaltingsChai} ; c'est un schéma sur~$\Spec(\Z[1/n])$. D'apr\`es~\cite[th. IV.5.7]{Deg@FaltingsChai}, il existe un sch\'ema semi-ab\'elien canonique sur~$\agnb$ qui prolonge le sch\'ema ab\'elien universel ; notons encore~$G$ ce sch\'ema semi-ab\'elien.
L'ensemble des points~$x\in\agnb\times\Spec(\F_p)$ tels que le groupe de $p$-torsion~$G_x[p]$ soit de rang multiplicatif~$p^r$ d\'efinit un sous-schéma localement ferm\'e
$$\agnb^r\:\hookrightarrow\: \agnb\times\Spec(\F_p)\: .$$
De m\^eme, l'ensemble des~$x\in\agnb\times\Spec(\F_p)$ tels que~$G_x[p]$ soit de rang multiplicatif~$\leq p^r$ d\'efinit un sous-schéma ferm\'e
$$\agnb^{\leq r}\:\hookrightarrow\:\agnb\times\Spec(\F_p)\: .$$
Nous verrons que~$\agnb^r$ est ouvert dans~$\agnb^{\leq r}$ ; il est dense d'apr\`es le th\'eor\`eme principal de~\cite{Newton@Oort}, mais nous n'utiliserons pas ce fait dans la suite.
Nous verrons qu'il existe un sous-groupe multiplicatif~$H^r_{\mathrm{mul}}$ de~$G[p]$ qui est fini, plat et totalement isotrope de rang~$p^r$ sur~$\agnb^r$ et qui est fibre \`a fibre le sous-groupe multiplicatif maximal de~$G[p]$.
Notons $$\agnb^{\mathrm{rig}}$$ la variété analytique rigide associ\'ee \`a~$\agnb\times\Spec(\Q_p)$, et notons
$]\agnb^r[$ et $]\agnb^{\leq r}[$ les tubes de~$\agnb^r$ et de $$\agnb^{\leq r}$$ dans~$\agnb^{\mathrm{rig}}$. D'apr\`es le th\'eor\`eme de rigidit\'e des tores~\sga{3}{iv}{th. 3.bis}, $H^r_{\mathrm{mul}}$ se rel\`eve canoniquement en un sous-groupe~$H^r_{\mathrm{can}}$ de~$G[p]$  qui est fini, plat et totalement isotrope de rang~$p^r$ sur~$]\agnb^r[$. Le but de cet article est de d\'emontrer que le sous-groupe~$H^r_{\mathrm{can}}$ surconverge, c'est-\`a-dire qu'il se prolonge \`a des voisinages stricts du tube. 

\begin{thprinc} Il existe un voisinage strict~$U$ de~$]\agnb^r[$ dans~$]\agnb^{\leq r}[$ 
et un sous-groupe~$H^r_{\mathrm{can},U}$ de~$G_U[p]$ qui est fini et plat de rang~$p^r$ sur~$U$ et qui prolonge~$H^r_\mathrm{can}$. Le sous-groupe~$H^r_{\mathrm{can},U}$ est unique et totalement isotrope d\`es que toute composante connexe de~$U$ rencontre~$]\agnb^r[$.
\end{thprinc}

On appelle~$H^r_{\mathrm{can},U}$ un \og sous-groupe canonique partiel de rang~$p^r$ sur~$U$\fg. Le sous-groupe canonique partiel de rang~$p^r$ est unique sur l'union des composantes connexes de~$U$ qui rencontrent~$]\agnb^r[$ ; remarquons d'ailleurs que cette union est encore un voisinage strict de~$]\agnb^r[$.

\begin{remarqueintro} Il n'existe pas de sous-groupe canonique partiel de rang~$p^r$ hors de~$]\agnb^{\leq r}[$. En effet, lorsque~$G[p]$ est de rang multiplicatif~$>p^r$, il n'admet aucun sous-groupe multiplicatif privil\'egi\'e de rang~$p^r$.
\end{remarqueintro}

Dans le cas o\`u~$r=g$, le sous-groupe canonique partiel n'est autre que le sous-groupe canonique habituel ; son existence a \'et\'e prouv\'ee par Lubin et Katz~\cite{Padic@Katz} lorsque~$g=1$ puis indépendamment par Abb\`es et Mokrane~\cite{Canonique@AbbesMokrane}, Andreatta et Gasbarri~\cite{Canonique@Andreatta} et Fargues~\cite{Harder@Fargues} pour~$g$ quelconque.
Ces auteurs définissent des filtrations sur le groupe de~$p$-torsion des schémas abéliens. Ils relient ensuite ces filtrations à l'invariant de Hodge des schémas abéliens. Cela leur permet d'obtenir l'équation de voisinages stricts où le sous-groupe canonique est défini.
Pour les variétés de Hilbert, Kisin et Lai~\cite{Canonique@KisinLai} ont construit le sous-groupe canonique par des méthodes de géométrie rigide. Notre travail s'inscrit dans cette lignée. Nos résultats sont utilisés par le premier auteur de cet article dans un travail en cours sur le prolongement analytique des formes de Siegel $p$-adiques : les sous-groupes canoniques partiels constituent un outil essentiel pour la compréhension de la dynamique des correspondances de Hecke.

Supposons \`a nouveau~$r$ quelconque et expliquons les grandes lignes de la d\'emonstration de l'assertion  d'existence du th\'eor\`eme. L'id\'ee g\'en\'erale est d'appliquer un th\'eor\`eme de surconvergence du tube d\^u \`a Berthelot~\cite{Rigide@Berthelot} \`a un morphisme d'oubli du niveau pour montrer que ce dernier r\'ealise un isomorphisme sur certains voisinages stricts ; le sous-groupe canonique partiel sera obtenu \textit{via} l'isomorphisme inverse.
Notons~$\agnd$ l'espace de modules des vari\'et\'es ab\'eliennes principalement polaris\'ees de genre~$g$ munies d'une structure de niveau principale en~$n$ et d'un sous-groupe fini, plat, totalement isotrope, tu\'e par~$p$ et de rang~$p^r$. Appelons~$H$ le sous-groupe fini et plat universel sur~$\agnd$. L'oubli de~$H$ induit un morphisme $\pi:\agnd\rightarrow\agn$. D'apr\`es le th\'eor\`eme principal de~\cite{Compact@Stroh}, il existe une compactification toro\"idale~$\agndb$ de~$\agnd$ associ\'ee au choix combinatoire ayant servi \`a d\'efinir~$\agnb$, et un morphisme~$\pib:\agndb\rightarrow \agnb$ \'etendant~$\pi$. Le groupe~$H$ se prolonge sur~$\agndb$ par adh\'erence sch\'ematique dans le groupe quasi-fini et plat~$G[p]$ ; nous obtenons ainsi un groupe quasi-fini et plat sur~$\agndb$ encore not\'e~$H$. D\'efinissons des sous-schémas de~$\agndb\times\Spec(\F_p)$ en posant
$$\begin{array}{lcllll}
\agndb^r & = & \{ \: s\in\agndb\times\Spec(\F_p) & | & H_{s}=G_{s}[p]^{\mathrm{mul}} & \}, \\
\agndb^{\leq r} & = & \{\: s\in\agndb\times\Spec(\F_p) & | & H_{s} \supset G_{s}[p]^{\mathrm{mul}} & \},\\
\agndb^{\geq r} & = & \{\: s\in\agndb\times\Spec(\F_p) & | & H_{s} \subset G_{s}[p]^{\mathrm{mul}} & \}.
\end{array}$$
Nous prouverons que~$H$ est fini et plat de rang~$p^r$ sur~$\agndb^{\leq r}$. Le morphisme~$\pib$ envoie~$\agndb^r$ sur~$\agnb^r$ et $$\agndb^{\leq r}$$ sur~$\agnb^{\leq r}$. Nous montrerons que~$\pib$ r\'ealise un isomorphisme entre~$\agndb^r$ et~$\agnb^r$ et est \'etale en tout point de $$\agndb^{\geq r}\: .$$
Nous pourrons donc appliquer le th\'eor\`eme~\cite[th. 1.3.5]{Rigide@Berthelot} au diagramme commutatif
$$\xymatrix{
\agndb^r\ar[d]^{\pib}_{\sim} \ar[r] &  \agndb^{\leq r} \ar[r] \ar[d]^{\pib} & \agndb^\mathrm{for} \ar[d]^{\pib} \\
\agnb^r \ar[r] & \agnb^{\leq r} \ar[r] & \agnb^\mathrm{for}}$$
dans lequel~$\agnb^{\mathrm{for}}$ et~$\agndb^{\mathrm{for}}$ d\'esignent les compl\'et\'es formels $p$-adiques de~$\agnb$ et de~$\agndb$ ; nous en d\'eduirons que~$\pib$ r\'ealise un isomorphisme
$$\pib\: : \:\: ]\agndb^r[\: \isolong \: ]\agnb^r[ $$
et qu'il existe des voisinages stricts~$U$ et~$V$ de~$]\agnb^r[$ et de~$]\agndb^r[$ dans~$]\agnb^{\leq r}[$ et dans~$]\agndb^{\leq r}[$ tels que~$\pib$ r\'ealise un isomorphisme de~$V$ sur~$U$. Il restera \`a vérifier que~$H$ est fini et plat de rang~$p^r$ sur~$V$ ; nous pourrons alors d\'efinir le sous-groupe canonique partiel de rang~$p^r$ sur~$U$ par la formule~$H^r_\mathrm{can}=(\pib^{-1})^*H$.

Expliquons \`a pr\'esent pourquoi~$\pib$ r\'ealise un isomorphisme de~$\agndb^r$ sur~$\agnb^r$ et est \'etale en tout point de $$\agndb^{\geq r}\: .$$
Ceci est \'evident sur les lieux modulaires obtenus par intersection avec~$\agnd$ et avec~$\agn$ : le caract\`ere \'etale r\'esulte par exemple du th\'eor\`eme de rigidit\'e des tores. Ainsi, il suffira de montrer que~$\pib$ est \'etale en tout point du bord de~$\agndb^{\geq r}$ ; cela occupera la majeure partie de l'article.
Nous montrerons en fait le caractère étale de~$\pib$ lorsque $n=1$ car cette hypothèse simplifie grandement la description explicite des compactifications toroïdales. Nous travaillerons alors non plus avec des schémas mais avec des champs algébriques. Le caractère étale de $\pib$ lorsque $n\geq 3$ suivra par changement de base.

Nous nous ram\`enerons \`a prouver un \'enonc\'e d'étalité sur les espaces de module de $1$-motifs. Nous utiliserons une description du morphisme d'oubli du niveau entre espaces de modules de $1$-motifs faisant intervenir un morphisme d'oubli du niveau entre vari\'et\'es de Siegel de genre~$< g$, une isog\'enie entre vari\'et\'es ab\'eliennes, et une isog\'enie entre tores. Le caract\`ere \'etale du morphisme d'oubli entre vari\'et\'es de Siegel de genre~$< g$ r\'esulte de la discussion pr\'ec\'edente, et nous montrerons que l'isog\'enie entre vari\'et\'es ab\'eliennes est \'etale, et que l'isog\'enie entre tores est un isomorphisme.

Remarquons que ni~$\agnb$, ni~$\agndb$, ni m\^eme le morphisme~$\agndb\rightarrow\agnb$ ne v\'erifient de propri\'et\'es modulaires simples. C'est la raison pour laquelle le caract\`ere \'etale est plus d\'elicat \`a d\'emontrer au bord : on ne peut pas appliquer directement de crit\`ere infinit\'esimal.

\section{Stratification des vari\'et\'es de Siegel par le rang multiplicatif}

Soit $g\geq 1$ un entier et~$p$ un nombre premier. Notons~$\ag$ le champ alg\'ebrique sur $\Spec(\Z)$ qui param\`etre les sch\'emas ab\'eliens principalement polaris\'es de genre~$g$, et notons~$G_0$ le sch\'ema ab\'elien universel (not\'e~$G$ dans l'introduction) sur~$\ag$. Soit $1\leq r \leq g$ un entier ; posons~$\mathscr{D}=\{r\}$. Notons~$\agd$ le champ sur~$\ag$ param\'etrant les sous-groupes~$H\subset G_0[p]$ qui sont finis, plats et totalement isotropes de rang~$p^r$. Notons~$H_0$ le sous-groupe universel (not\'e~$H$ dans l'introduction) sur~$\agd$ et $G_1=G_0/H_0$. D'apr\`es~\cite{Gamma0@DeJong}, la donn\'ee du sch\'ema ab\'elien principalement polaris\'e~$G_0$ et du sous-groupe totalement isotrope $H_0\subset G_0[p]$ est \'equivalente \`a celle d'une cha\^ine d'isog\'enie
$$G_0 \: \stackrel{\alpha}{\longrightarrow} \: G_1 \: \stackrel{\lambda}{\longrightarrow} \: G_1^t \: \stackrel{\alpha^t}{\longrightarrow} \: G_0^t \: ,$$
o\`u~$\alpha^t$ est duale de~$\alpha$, o\`u~$\lambda$ est une polarisation de~$G_1$, et o\`u la compos\'ee~$\alpha^t \circ \lambda \circ \alpha$ est \'egale \`a la multiplication par~$p$ de la polarisation principale de~$G_0$. Nous pr\'ef\`ererons fr\'equemment parler de cha\^ines d'isog\'enies plut\^ot que de sous-groupes : les conditions d'isotropie s'expriment de mani\`ere plus synth\'etique dans le langage des cha\^ines.

Soit $S$ le spectre d'un corps de caract\'eristique~$p$ et~$X$ un groupe fini de $p$-torsion sur~$S$. 
Il existe un d\'evissage canonique
$$0 \:\longrightarrow\: X^{\mathrm{mul}} \: \longrightarrow \: X \:\longrightarrow\: X^{\mathrm{non-mul}} \:\longrightarrow\: 0$$
en groupes finis sur~$S$ tels que $X^{\mathrm{mul}}$ soit multiplicatif et que $X^{\mathrm{non-mul}}$ n'admette pas de sous-groupe multiplicatif. On appelle rang multiplicatif de~$X$ le rang de~$X^{\mathrm{mul}}$. 

Dans le cas o\`u~$S$ est un sch\'ema sur~$\F_p$ qui n'est pas n\'ecessairement le spectre d'un corps, et~$X$ est un groupe fini et plat de $p$-torsion sur~$S$, la fonction de~$S$ dans~$\mathbb{N}$ qui \`a $s\in S$ associe le rang multiplicatif de $X_s$ est semi-continue sup\'erieurement :
elle peut d\'ecro\^itre lors de sp\'ecialisations.

D\'efinissons plusieurs sous-sch\'emas localement ferm\'es r\'eduits de~$\ag\times\Spec(\F_p)$ et de \linebreak $\agd\times\Spec(\F_p)$ de la mani\`ere suivante :
$$\begin{array}{lcllll}
\ag^r & = & \{ \: s\in \ag\times\Spec(\F_p) & |&  \mathrm{rg}(G_{0,s}[p]^{\mathrm{mul}})=p^r & \}, \\
\ag^{\geq r} & = & \{ \: s\in\ag\times\Spec(\F_p) & |&  \mathrm{rg}(G_{0,s}[p]^{\mathrm{mul}})\geq p^r & \}, \\
\ag^{\leq r} & = & \{ \: s\in\ag\times\Spec(\F_p) & | & \mathrm{rg}(G_{0,s}[p]^{\mathrm{mul}})\leq p^r & \}, \\
\agd^r & = & \{ \: s\in\agd\times\Spec(\F_p) & | & H_{0,s}=G_{0,s}[p]^{\mathrm{mul}} & \}, \\
\agd^{\geq r} & = & \{\: s\in\agd\times\Spec(\F_p) &  | & H_{0,s} \subset G_{0,s}[p]^{\mathrm{mul}} & \}, \\
\agd^{\leq r} & = & \{\: s\in\agd\times\Spec(\F_p) & | & H_{0,s} \supset G_{0,s}[p]^{\mathrm{mul}} & \}.
\end{array}$$
Remarquons que $\pi:\agd\rightarrow\ag$ envoie~$\agd^r$ dans~$\ag^r$, $\agd^{\geq r}$ dans~$\ag^{\geq r}$ et $\agd^{\leq r}$ dans~$\ag^{\leq r}$.
On d\'eduit de la semi-continuit\'e du rang multiplicatif que les immersions $$\ag^r\hookrightarrow \ag^{\geq r}$$ et $$\agd^r\hookrightarrow\nolinebreak \agd^{\geq r}$$ sont ferm\'ees et que les immersions $\ag^r\hookrightarrow \ag^{\leq r}$ et $\agd^r\hookrightarrow \agd^{\leq r}$ sont ouvertes. De m\^eme, les immersions $$\ag^{\geq r} \hookrightarrow \ag\times\Spec(\F_p)$$ et $$\agd^{\geq r} \hookrightarrow \agd\times\Spec(\F_p)$$ sont ouvertes et les immersions $$\ag^{\leq r} \hookrightarrow \ag\times\Spec(\F_p)$$ et $$\agd^{\leq r} \hookrightarrow \agd\times\Spec(\F_p)$$ sont ferm\'ees. D'apr\`es le th\'eor\`eme de Oort~\cite{Newton@Oort},~$\ag^r$ est dense dans~$\ag^{\leq r}$. Nous n'utiliserons pas ce fait dans la suite. 

\begin{exemple} Si $r=g$, alors $\ag^r=\ag^{\geq r}$ est le lieu ordinaire de~$\ag\times \Spec(\F_p)=\ag^{\leq r}$. \linebreak De m\^eme, $$\agd^r=\agd^{\geq r}$$ est le lieu ordinaire-multiplicatif et $\agd^{\leq r}$ est la composante irr\'eductibles de~$\agd\times\Spec(\F_p)$ o\`u le sous-groupe universel~$H_0$ est \'egal au noyau du morphisme de Frobenius de~$G_0$.
\end{exemple}

Rappelons qu'un hom\'eomorphisme universel est un morphisme de sch\'emas qui induit une homéomorphisme apr\`es tout changement de base.

\begin{lemme}\label{lemme_homeo_univ} La projection canonique $\pi:\agd^r\rightarrow\ag^r$ est un hom\'eomorphisme universel.
\end{lemme}

\begin{demo} On vérifie aisément que $\pi$ est propre. Il suffit donc de montrer qu'elle est universellement bijective.
Soit $S$ un sch\'ema sur~$\ag^r$. Si l'on note $$S'=S\times_{\ag^r}\agd^r\: ,$$ on obtient un sch\'ema ab\'elien~$G$ sur~$S$ et un sous-groupe~$H$ de~$G_{S'}[p]$. Pour tout~$s\in S$, le rang multiplicatif de~$G_s[p]$ est~$p^r$. Le groupe~$G_s[p]$ contient donc un unique sous-groupe multiplicatif de rang~$p^r$, \`a savoir~$G_s[p]^{\mathrm{mul}}$. Comme l'ensemble sous-jacent \`a $\pi^{-1}(s)$ est en bijection avec celui des sous-groupes multiplicatifs de rang~$p^r$ de~$G_s[p]$, on voit que~$\pi$ induit une bijection entre~$S$ et~$S'$.\end{demo}

\begin{proposition} \label{prop_pi_etale_ouvert} La projection canonique $\pi:\agd\rightarrow \ag$ est \'etale en tout point de~$\agd^{\geq r}$ et induit une surjection de~$\agd^{\geq r}$ sur~$\ag^{\geq r}$.
\end{proposition}

\begin{demo} La surjectivit\'e est \'evidente. Pour montre l'\'etalit\'e, on applique le crit\`ere infinit\'esimal.
Soit~$S$ un sch\'ema affine et~$S'$ un sous-sch\'ema ferm\'e de~$S$ d\'efini par un id\'eal de carr\'e nul. On se donne un diagramme commutatif
$$\xymatrix{S' \ar[r]\ar[d] & \agd^{\geq r} \ar[d] \\
S \ar[r] & \ag}$$
et l'on veut montrer qu'il existe une unique factorisation $S\rightarrow \agd$. Le diagramme pr\'ec\'edent induit un sch\'ema ab\'elien~$G$ sur~$S$ et un groupe fini et plat~$H\subset G_{S'}[p]$ dont toutes les fibres sont multiplicatives. Le sous-groupe~$H$ de~$G_{S'}$ est donc multiplicatif sur~$S'$, et s'\'etend de mani\`ere unique en un sous-groupe multiplicatif de~$G$ par le th\'eor\`eme de rigidit\'e des tores~\sga{3}{iv}{th. 3.bis}. On a bien obtenu une factorisation $S\rightarrow \agd$.
\end{demo}

Le corollaire suivant r\'esulte du lemme~\ref{lemme_homeo_univ} et de la proposition~\ref{prop_pi_etale_ouvert} en remarquant que tout hom\'eomorphisme universel \'etale est un isomorphisme.

\begin{corollaire} Le morphisme d'oubli du niveau induit un isomorphisme $\pi:\agd^r\iso\ag^r$.
\end{corollaire}
 
L'isomorphisme inverse $\ag^r \rightarrow \agd^r$ permet d'obtenir un sous-groupe $H_0\subset G_0[p]\times_{\ag} \ag^r$ fini, plat et multiplicatif, tel que pour tout $s\in\ag^r$, on ait $$H_{0,s}\: = \: G_{0,s}[p]^{\mathrm{mul}}\: .$$
Nous aurions pu d\'eduire cela directement de~\cite[coro.~II.1.2]{Shimura@HarrisTaylor}.

\begin{remarque} Dans le cas o\`u $r=g$, le groupe~$H_0$ co\"incide avec le noyau du morphisme de Frobenius de~$G_0\times\ag^r$.
\end{remarque}

\section{Stratification sur les compactifications toro\"idales}

\subsection{Rappels} Dans ce paragraphe, nous rappelons les propri\'et\'es des compactifications toro\"idales qui seront utilis\'ees dans la suite. Commençons par fixer quelques notations. 

Soit $V=\oplus_{i=1}^{2g}\: \Z \: e_i$ un $\Z$-module libre de rang~$2g$ muni d'une base~$(e_i)$. Munissons~$V$ de l'accouplement symplectique donn\'e par la formule
$$\langle \: \sum_{i=1}^{2g} a_i \: e_i \: , \: \sum_{j=1}^{2g} b_j \: e_j \: \rangle \: = \: a_1 b_{2g} - a_{2g} b_1 + a_2 b_{2g-1} - a_{2g-1} b_2 + \cdots + a_{g} b_{g+1} - a_{g+1} b_g \: .$$
Notons $\mathfrak{C}$ l'ensemble des sous-modules facteurs directs totalement isotropes de~$V$. Pour tout~$V'\in\mathfrak{C}$, notons $C(V/V'^\perp)$ le c\^one des formes quadratiques semi-d\'efinies positives \`a radical rationnel sur $V/V'^\perp\otimes \R$. Toute inclusion $V''\subset V'$ entre \'el\'ements de~$\mathfrak{C}$ induit une inclusion $C(V/V''^\perp)\subset C(V/V'^\perp)$. D\'efinissons un complexe conique en posant
$$\mathscr{C} \: = \: \left( \coprod_{V'\in \mathfrak{C}} \: C(V/V'^\perp) \right) / \sim$$
o\`u $\sim$ est la relation d'\'equivalence engendr\'ee par les inclusions $C(V/{V''}^\perp )\subset C(V/V'^\perp )$ pour tous $V'', \: V'\in\mathfrak{C}$ v\'erifiant~$V''\subset V'$.

Notons $\Gamma=\gsp(\Z)$ le groupe des similitudes symplectiques de~$V$. Ce groupe agit d'une mani\`ere compatible sur~$\mathfrak{C}$ et sur~$\mathscr{C}$.
Associons \`a~$\mathscr{D}=\{r\}$ une cha\^ine de~$V$ en posant
\begin{eqnarray*}
\mathcal{V}_\mathscr{D}^\bullet & = & (\mathcal{V}_\mathscr{D}^0\subset \mathcal{V}_\mathscr{D}^1 \subset \mathcal{V}_\mathscr{D}^2 \subset \mathcal{V}_\mathscr{D}^3) \\
\mathcal{V}_\mathscr{D}^0 & = & p\cdot V \\
\mathcal{V}_\mathscr{D}^1 & = & ( \oplus_{i=1}^{r}\: \Z \cdot e_i) \: \oplus \: (\oplus_{i>r}\:\: p\:\Z \cdot e_i )\\
\mathcal{V}_\mathscr{D}^2 & = & \{ \: v \in V \: | \: <v,\mathcal{V}_\mathscr{D}^1>\: \in\: p\:\Z \: \} \\
\mathcal{V}_\mathscr{D}^3 & = & V
\end{eqnarray*}
Le stabilisateur~$\Gamma_\mathscr{D}$ de $\mathcal{V}_\mathscr{D}^\bullet$ dans~$\Gamma$ est le groupe parahorique habituellement associ\'e \`a~$\mathscr{D}$ : il est form\'e des matrices de~$\Gamma$ dont la r\'eduction modulo~$p$ est triangulaire sup\'erieure par blocs, le premier bloc \'etant de taille~$r\times r$.

Choisissons d\'esormais une d\'ecomposition poly\'edrale rationnelle $\Gamma$-admissible~$\mathscr{S}$ de~$\mathscr{C}$ comme dans~\cite[par. 3.2.3]{Compact@Stroh}. Cette d\'ecomposition est en particulier~$\Gamma_\mathscr{D}$-admissible. Les constructions effectu\'ees dans~\cite[ch. IV]{Deg@FaltingsChai} et dans~\cite{Compact@Stroh} permettent d'associer \`a~$\mathscr{S}$ deux champs alg\'ebriques~$\agb$ et~$\agdb$ propres sur~$\Spec(\Z)$, ainsi qu'un diagramme cart\'esien
\begin{equation}
\xymatrix{\agd \ar[r] \ar[d]^{\pi} & \agdb \ar[d]^{\pib} \\
\ag \ar[r] & \agb }
\end{equation}\label{diagramme_ag}
dans lequel les fl\`eches horizontales sont des immersions ouvertes d'image dense. Les champs~$\agb$ et~$\agdb$ sont munis de stratifications respectivement param\'etr\'ees par~$\mathfrak{C}/\Gamma$ et par~$\mathfrak{C}/\Gamma_\mathscr{D}$. Le morphisme~$\pib$ respecte ces stratifications, et la fl\`eche induite par~$\pib$ sur l'ensemble des strates est la projection canonique~$\mathfrak{C}/\Gamma\rightarrow\mathfrak{C}/\Gamma_\mathscr{D}$.

\begin{remarque} Les champs~$\agb$ et~$\agdb$ sont \'egalement munis de stratifications param\'etr\'ees par les ensembles finis~$\mathscr{S}/\Gamma$ et~$\mathscr{S}/\Gamma_\mathscr{D}$. Ces stratifications raffinent celles param\'etr\'ees par~$\mathfrak{C}/\Gamma$ et~$\mathfrak{C}/\Gamma_\mathscr{D}$ d'une mani\`ere compatible au morphisme de~$\mathscr{S}$ dans~$\mathfrak{C}$ qui envoie~$\sigma$ sur l'unique~$V'$ tel que~$\sigma$ soit inclus dans l'int\'erieur de~$C(V/V'^\perp)$. Dans la suite de l'article, nous n'utiliserons pas les stratifications d'origine combinatoire param\'etr\'ees par~$\mathscr{S}/\Gamma$ et~$\mathscr{S}/\Gamma_\mathscr{D}$.
\end{remarque}

\subsubsection*{Description du bord des compactifications} Rappelons que l'on peut donner un sens au hens\'elis\'e d'un champ alg\'ebrique~$X$ le long d'un sous-champ ferm\'e~$Y$ : c'est la limite projective des champs alg\'ebriques munis d'un morphisme affine \'etale vers~$X$ induisant un isomorphisme au-dessus de~$Y$. Cette limite projective existe dans la $2$-cat\'egorie des champs alg\'ebriques puisque les fl\`eches de transition sont affines. Dans le cas o\`u~$X$ est un sch\'ema et~$Y$ un point ferm\'e, la d\'efinition pr\'ec\'edente co\"incide avec celle du hens\'elis\'e habituel, car les ouverts \'etales affines sont cofinaux parmis tous les ouverts \'etales. Dans le cas g\'en\'eral, la compl\'etion formelle de~$X$ le long de~$Y$ co\"incide avec celle du hens\'elis\'e de~$X$ le long de~$Y$. Notons $X^{h_Y}$ le hens\'elis\'e de~$X$ le long de~$Y$.

Soit~$V'\in\mathfrak{C}$. Notons~$Z_{V'}$ et~$Z_{V',\mathscr{D}}$ les strates de~$\agb$ et de~$\agdb$ param\'etr\'ees par~$V'$. Nous appelons \og hens\'elis\'e du diagramme~(\ref{diagramme_ag}) le long de la~$V'$-strate \fg le diagramme commutatif
$$\xymatrix{\agd \times_{\agdb} \agdb^{h_{Z_{V',\mathscr{D}}}} \ar[r] \ar[d]^{\pi} & \agdb^{h_{Z_{V',\mathscr{D}}}} \ar[d]^{\pib} \\
\ag \times_{\agb} \agb^{h_{Z_{V'}}} \ar[r] & \agb^{h_{Z_{V'}}} }$$
De m\^eme pour tout diagramme cart\'esien stratifi\'e analogue \`a (\ref{diagramme_ag}).
D'apr\`es~\cite[th. IV.5.7]{Deg@FaltingsChai} et~\cite[th. 3.2.7.1]{Compact@Stroh}, le hens\'elis\'e du diagramme~(\ref{diagramme_ag}) le long de la~$V'$-strate est canoniquement isomorphe au quotient par un groupe discret du hens\'elis\'e d'un diagramme
\begin{equation} \label{diagramme_m}
\xymatrix{\mvd \ar[r] \ar[d]^{\pi_{V'}} & \mvdb \ar[d]^{\pib_{V'}} \\
\mv \ar[r] & \mvb }
\end{equation}
le long d'une strate ; toutes les notations employ\'ees seront expliqu\'ees dans la suite de ce paragraphe.
On en d\'eduit une description du compl\'et\'e formel de~$\agb$ et de~$\agdb$ le long des~strates.

Il nous faut \`a pr\'esent expliquer le diagramme~(\ref{diagramme_m}). En un mot : $\mv$ est l'espace de modules des $1$-motifs de genre $g$ principalement polaris\'es et rigidifi\'es par~$V'$~\cite[d\'ef. 1.4.3]{Compact@Stroh}, et le morphisme $\mvd\rightarrow\mv$ param\`etre les structures de niveau parahoriques sur le~$1$-motif universel qui sont de type~$\mathscr{D}$ en~$p$ et rigidifi\'ees par~$V'$~\cite[par. 1.4.9]{Compact@Stroh}. Le champ~$\mvd$ param\`etre donc certains sous-groupes finis, plats et totalement de rang~$p^r$ du groupe des points de~$p$-torsion du $1$-motif universel~\cite[par. 1.2.4]{Compact@Stroh}. Rappelons que le groupe des points de~$p$-torsion d'un $1$-motif est muni d'une filtration \`a trois crans dont les gradu\'es sont appel\'es \og partie  torique \fg, \og partie ab\'elienne \fg et \og partie p\'eriodique \fg. Le champ~$\mvd$ param\`etre les sous-groupes~$H$ dont les parties toriques, ab\'eliennes et p\'eriodiques sont d\'etermin\'ees par la position relative de la cha\^ine~$\mathcal{V}_\mathscr{D}^\bullet$ et du drapeau~$(V'\subset V'^\perp \subset V)$ : l'existence d'une rigidification par~$V'$ impose notamment que le rang torique de~$H$ soit \'egal au cardinal du noyau de
$$(\mathcal{V}_\mathscr{D}^0\cap V')\otimes \F_p \:\longrightarrow\: (\mathcal{V}_\mathscr{D}^1\cap V')\otimes \F_p$$
(qui est \'egal au cardinal du noyau de $(\mathcal{V}_\mathscr{D}^2 / \:\mathcal{V}_\mathscr{D}^2\cap V'^\perp)\otimes \F_p \:\longrightarrow\: (\mathcal{V}_\mathscr{D}^3 /\:\mathcal{V}_\mathscr{D}^3\cap V'^\perp)\otimes \F_p$ comme on le voit par dualit\'e), 
que le rang ab\'elien de~$H$ soit \'egal au cardinal du noyau de
$$(\mathcal{V}_\mathscr{D}^0\cap V'^\perp /\: \mathcal{V}_\mathscr{D}^0\cap V')\otimes \F_p \:\longrightarrow\: (\mathcal{V}_\mathscr{D}^1\cap V'^\perp/\:\mathcal{V}_\mathscr{D}^1\cap V')\otimes \F_p\: ,$$
et que le rang p\'eriodique de~$H$ soit \'egal au cardinal du noyau de
$$(\mathcal{V}_\mathscr{D}^0 / \:\mathcal{V}_\mathscr{D}^0\cap V'^\perp)\otimes \F_p \:\longrightarrow\: (\mathcal{V}_\mathscr{D}^1 /\:\mathcal{V}_\mathscr{D}^1\cap V'^\perp)\otimes \F_p\: .$$
Terminons l'explication rapide le diagramme~(\ref{diagramme_m}) : les immersions ouvertes $\mv\hookrightarrow \mvb$ et $\mvd\hookrightarrow\nolinebreak\mvdb$ sont des plongements toro\"idaux localement de type fini associ\'es \`a la restriction de~$\mathscr{S}$ \`a~$C(V/V'^\perp)$.

Passons \`a une description plus pr\'ecise des champs alg\'ebriques~$\mv$ et~$\mvd$ : il existe un diagramme commutatif de champs alg\'ebriques sur~$\Spec(\Z)$
\begin{equation}\label{diagramme_mba}
\xymatrix{
\mvd \ar[r] \ar[d]^{\pi_{V'}} & \bvd \ar[d]\ar[r] & \avd \ar[d]\\
\mv \ar[r] & \bv \ar[r] & \av}
\end{equation}
dans lequel
\begin{enumerate}
\item $\av$ est la vari\'et\'e de Siegel qui param\`etre les sch\'emas ab\'eliens principalement polaris\'es de genre~$g-\mathrm{rg}(V')$,
\item $\avd$ est la vari\'et\'e de Siegel qui param\`etre les sch\'emas ab\'eliens principalement polaris\'es de genre~$g-\mathrm{rg}(V')$ munis d'un sous-groupe fini, plat et totalement isotrope de rang \'egal au cardinal du noyau de
$$(\mathcal{V}_\mathscr{D}^0\cap V'^\perp /\: \mathcal{V}_\mathscr{D}^0\cap V')\otimes \F_p \:\longrightarrow\: (\mathcal{V}_\mathscr{D}^1\cap V'^\perp/\:\mathcal{V}_\mathscr{D}^1\cap V')\otimes \F_p\: ,$$
\item $\bv\rightarrow\av$ param\`etre les extensions
$$ 0 \: \longrightarrow T_0 \:\longrightarrow\: \tilde{G}_0 \:\longrightarrow\: A_0 \:\longrightarrow\: 0 \: , $$
o\`u $T_0=\mathrm{Hom}(V/V'^\perp,\Gm)$ et~$A_0$ est la vari\'et\'e ab\'elienne universelle sur~$\av$,
\item $\bvd\rightarrow\avd$ param\`etre les diagrammes commutatifs d'extensions
$$\xymatrix{
0 \ar[r] &  T_0 \ar[r]\ar[d] & \tilde{G}_0 \ar[r]\ar[d] & A_0 \ar[r] \ar[d] & 0 \\
0 \ar[r] &  T_1 \ar[r]\ar[d] & \tilde{G}_1 \ar[r]\ar[d] & A_1 \ar[r] \ar[d] & 0 \\
0 \ar[r] &  T'_1 \ar[r]\ar[d] & \tilde{G}'_1 \ar[r]\ar[d] & A_1^t \ar[r] \ar[d] & 0 \\
0 \ar[r] &  T'_0 \ar[r] & \tilde{G}'_0 \ar[r] & A_0^t \ar[r] & 0}
$$
o\`u la cha\^ine $T_0\rightarrow T_1\rightarrow T'_1 \rightarrow T'_0$ est fix\'ee par les formules
$T_i=\mathrm{Hom}(\mathcal{V}_\mathscr{D}^{3-i}/V'^\perp,\Gm)$ et  
$T'_i=\mathrm{Hom}(\mathcal{V}_\mathscr{D}^{i}/V'^\perp,\Gm)$ pour~$i=0$ et $1$, o\`u
la cha\^ine $A_0\rightarrow A_1\rightarrow A^t_1 \rightarrow A^t_0$ est fix\'ee et provient de la cha\^ine d'isog\'enies universelles sur~$\avd$, et o\`u le compos\'e vertical $\tilde{G}_0\rightarrow\tilde{G}'_0$ est la mutiplication d'un isomorphisme par~$p$ et $A_0\rightarrow A_0^t$ est la multiplication de la polarisation principale par~$p$,
\item $\mv\rightarrow\bv$ param\`etre les $1$-motifs principalement polaris\'es $[V/V'^\perp\rightarrow \tilde{G}_0]$,
\item $\mvd\rightarrow\bvd$ param\`etre les diagrammes commutatifs de $1$-motifs
$$\xymatrix{
[ \: \mathcal{V}_\mathscr{D}^0 /  V'^\perp \ar[d]\ar[r] & \tilde{G}_0 \ar[d] \: ] \\
[ \: \mathcal{V}_\mathscr{D}^1 /  V'^\perp \ar[d]\ar[r] & \tilde{G}_1 \ar[d] \: ] \\
[ \: \mathcal{V}_\mathscr{D}^2 /  V'^\perp \ar[d]\ar[r] & \tilde{G}'_1 \ar[d] \: ] \\
[ \: \mathcal{V}_\mathscr{D}^3 /  V'^\perp \ar[r] & \tilde{G}'_0  \: ]}
$$
o\`u la cha\^ine $\mathcal{V}_\mathscr{D}^0 /  V'^\perp$ est fix\'ee et est d\'etermin\'ee par~$\mathcal{V}_\mathscr{D}^\bullet$, o\`u la cha\^ine $\tilde{G}_0 \rightarrow \tilde{G}_1 \rightarrow \tilde{G}'_1 \rightarrow \tilde{G}'_0$ est fix\'ee et provient de la cha\^ine d'isog\'enies universelles sur~$\bvd$, o\`u le $1$-motif $[\mathcal{V}_\mathscr{D}^3 /  V'^\perp \rightarrow \tilde{G}'_0]$ est dual de $[\mathcal{V}_\mathscr{D}^0 /  V'^\perp \rightarrow \tilde{G}_0]$ et la compos\'ee des fl\`eches verticales est la mutiplication d'une polarisation principale de $[\mathcal{V}_\mathscr{D}^0 /V'^\perp \rightarrow \tilde{G}_0]$ par~$p$.
\end{enumerate}

Dans le diagramme~(\ref{diagramme_mba}), les fl\`eches verticales correspondent \`a l'oubli du niveau, et sont d\'efinies de mani\`ere \'evidente en utilisant l'isomorphisme~$\mathcal{V}_\mathscr{D}^0\iso V$ de multiplication par~$p$. Le morphisme $\bv\rightarrow \av$ est un sch\'ema ab\'elien par la formule de Barsotti-Weil, et le morphisme $\mv\rightarrow\bv$ est un torseur sous le tore~$\Hom(\mathrm{Sym}^2(V/V'^\perp),\Gm)$. D'apr\`es~\cite[par. 1.4.7]{Compact@Stroh}, le morphisme $\bvd\rightarrow\agd$ est \'egalement un sch\'ema ab\'elien et $\mvd\rightarrow\bvd$ un torseur sous un tore.

Enfin, $\mv\hookrightarrow\mvb$ et $\mvd\hookrightarrow\mvdb$ sont les fibr\'es en plongements toro\"idaux sur~$\bv$ et~$\bvd$ associ\'es \`a la restriction de~$\mathscr{S}$ \`a~$C(V/V'^\perp)$. Ce sont des champs alg\'ebriques localement de type fini sur~$\Spec(\Z)$ munis de stratifications respectivement param\'etr\'ees par $\{V''\in\mathfrak{C}\: | \: V''\subset\nolinebreak V'\}/\Gamma$ et  par $\{V''\in\mathfrak{C}\: | \: V''\subset V'\}/\Gamma_\mathscr{D}$.
Notons~$Y_{V'}$ et~$Y_{V',\mathscr{D}}$ les strates param\'etr\'ees par~$V'$ : ce sont les uniques strates ferm\'ees.
Notons~$\Gamma_{V'}$ le stabilisateur de~$V'$ dans~$\Gamma$, et $\Gamma_{V',\mathscr{D}}=\Gamma_{V'}\cap \Gamma_\mathscr{D}$. Le groupe discret~$\Gamma_{V'}$ agit sans point fixe sur~$\mvb$ en pr\'eservant~$Y_{V'}$, et le quotient~$\mv/\Gamma_{V'}$ est de type fini sur~$\Spec(\Z)$. Il en est de m\^eme pour~$\Gamma_{V',\mathscr{D}}$ et~$\mvdb$. Nous pouvons \`a pr\'esent pr\'eciser le lien entre compactifications toro\"idales et espaces de modules de $1$-motifs. La proposition suivante r\'esulte de~\cite[th. IV.5.7]{Deg@FaltingsChai} et de~\cite[th. 3.2.7.1]{Compact@Stroh}.

\begin{proposition} \label{prop_struc_locale} Il existe des isomorphismes canoniques stratifi\'es
$$\mvb^{h_{Y_{V'}}} \: / \: \Gamma_{V'} \: \isolong \: \agb^{h_{Z_{V'}}} \quad \mathrm{et} \quad \mvdb^{h_{Y_{V',\mathscr{D}}}} \: / \: \Gamma_{V',\mathscr{D}} \: \isolong \: \agdb^{h_{Z_{V',\mathscr{D}}}}\: .$$
En particulier, il existe des isomorphismes $Y_{V'}/\Gamma_{V'}\iso Z_{V'}$ et $Y_{V',\mathscr{D}}/\Gamma_{V',\mathscr{D}}\iso Z_{V',\mathscr{D}}$.
\end{proposition}

\subsection{D\'efinition de la stratification} Nous prolongeons \`a pr\'esent aux compactifications toro\"idales les sous-sch\'emas d\'efinis dans la premi\`ere partie. Nous raisonnerons pour cela strate de bord par strate de bord, en expliquant quelles strates et quels sous-sch\'emas de chacune de ces strates consid\'erer. Nous verrons \`a la fin de ce paragraphe que les sous-sch\'emas construits peuvent \^etre d\'ecrit comme le lieu des compactifications toro\"idales o\`u le rang multiplicatif d'un groupe quasi-fini et plat est fix\'e ; cela fera le lien avec les d\'efinitions de l'introduction de l'article.

Introduisons au pr\'ealable quelques notations :
$$\begin{array}{lll}
\mathfrak{C}^r & = & \{ \: V'\in\mathfrak{C} \:\:\: |\: \mathrm{rg}(V')\leq r \: \} \\
\mathfrak{C}^{\geq r} & = & \mathfrak{C}  \\
\mathfrak{C}^{\leq r} & = & \mathfrak{C}^r  \\
\mathfrak{C}_\mathscr{D}^r & = & \{\: V'\in\mathfrak{C}^r \: |\:  \mathcal{V}_\mathscr{D}^0/V'^\perp=\mathcal{V}_\mathscr{D}^1/V'^\perp=\mathcal{V}_\mathscr{D}^2/V'^\perp \:\}\\
\mathfrak{C}_\mathscr{D}^{\geq r} & = & \{ \: V'\in\mathfrak{C} \:\:\: | \: \mathcal{V}_\mathscr{D}^0/V'^\perp=\mathcal{V}_\mathscr{D}^1/V'^\perp \: \} \\
\mathfrak{C}_\mathscr{D}^{\leq r} & = & \mathfrak{C}_\mathscr{D}^r \: .
\end{array}$$
Pour tout $V'\in\mathfrak{C}^r$ (resp. $\mathfrak{C}^{\geq r}$, $\mathfrak{C}^{\leq r}$), appelons $\av^r$ (resp. $\av^{\geq r}$, $\av^{\leq r}$) le sous-sch\'ema localement ferm\'e de~$\av\times\Spec(\F_p)$ o\`u le rang multiplicatif est \'egal (resp. sup\'erieur, inf\'erieur) \`a~$p^{r-\mathrm{rg}(V')}$. Avec les notations de la premi\`ere partie, on a donc 
$$\av^r=\mathcal{A}_{g-\mathrm{rg}(V')}^{r-\mathrm{rg}(V')}\quad \quad \av^{\geq r}=\mathcal{A}_{g-\mathrm{rg}(V')}^{\geq r-\mathrm{rg}(V')}\quad \quad  \av^{\leq r}=\mathcal{A}_{g-\mathrm{rg}(V')}^{\leq r-\mathrm{rg}(V')} \: .$$
De m\^eme, notons~$\avd^r$ le sous-sch\'ema localement ferm\'e de~$\avd\times\Spec(\F_p)$ o\`u le rang multiplicatif est~$p^{r-\nolinebreak\mathrm{rg}(V')}$ et o\`u~$\mathrm{Ker}(A_0 \rightarrow A_1)$ est multiplicatif, et ainsi de suite pour~$\avd^{\geq r}$ et~$\avd^{\leq r}$. 
D\'efinissons $\bv^r$, $\bv^{\geq r}$, $\bv^{\leq r}$, $\bvd^r$, ..., $Y_{V'}^r$, ..., $Y_{V',\mathscr{D}}^{r}$, ..., $\mv^r$, ..., $\mvdb^r$, ... par produits fibr\'es. Par exemple :
$$
\begin{array}{ccc}
\bv^r=\bv\times_{\av} \av^r & \bv^{\geq r}=\bv\times_{\av} \av^{\geq r} & \bvd^{r}=\bvd\times_{\avd} \avd^{r} \\
Y_{V'}^r=Y_{V'}\times_{\av} \av^r & Y_{V',\mathscr{D}}^{r}=Y_{V',\mathscr{D}}\times_{\avd} \avd^{r} & \mvdb^r=\mvdb\times_{\avd}\avd^r\: .
\end{array}$$
Les immersions $\mvb^r\hookrightarrow \mvb^{\geq r}$, $\mvdb^r\hookrightarrow \mvdb^{\geq r}$, $\mvb^{\leq r}\hookrightarrow \mvb \times \Spec(\F_p)$ et $\mvdb^{\leq r}\hookrightarrow \mvdb \times \Spec(\F_p)$ sont ferm\'ees ; les immersions $\mvb^r\hookrightarrow \mvb^{\leq r}$, $\mvdb^r\hookrightarrow \mvdb^{\leq r}$, $\mvb^{\geq r}\hookrightarrow \mvb \times \Spec(\F_p)$ et $\mvdb^{\geq r}\hookrightarrow \mvdb \times \Spec(\F_p)$ sont ouvertes.

\begin{remarque} En vertu du théorème de Oort~\cite{Newton@Oort}, on peut voir que $\mvb^r$ est dense dans~$\mvb^{\leq r}$ ; nous n'utiliserons pas ce fait dans la suite.
\end{remarque}

En quotientant~$Y_{V'}^r$, $Y_{V'}^{\geq r}$ et~$Y_{V'}^{\leq r}$ par~$\Gamma_{V'}$ (resp. $Y_{V',\mathscr{D}}^r$, $Y_{V',\mathscr{D}}^{\geq r}$ et~$Y_{V',\mathscr{D}}^{\leq r}$ par~$\Gamma_{V',\mathscr{D}}$), l'on obtient des sous-sch\'emas localement ferm\'es~$Z_{V'}^r$, $Z_{V'}^{\geq r}$ et~$Z_{V'}^{\leq r}$ de~$\agb\times \Spec(\F_p)$ (resp. $Z_{V',\mathscr{D}}^r$, $Z_{V',\mathscr{D}}^{\geq r}$ et~$Z_{V',\mathscr{D}}^{\leq r}$ de~$\agdb\times \Spec(\F_p)$.

\begin{remarque} Lorsque~$V'=0$, on a $Z_{V'}=\ag$, $Z_{V'}^r=\ag^r$, $Z_{V',\mathscr{D}}=\agd$, $Z_{V',\mathscr{D}}^r=\agd^r$ et ainsi de suite.
\end{remarque}

Nous prolongeons les sous-sch\'emas localement ferm\'es d\'efinis dans la premi\`ere partie aux compactifications toro\"idales en posant
$$\begin{array}{cccccccccccccc}
\agb^r & = & \bigcup_{V'\in\mathfrak{C}^r} & Z_{V'}^r & \:\:  & \agb^{\geq r} & = & \bigcup_{V'\in\mathfrak{C}^{\geq r}} & Z_{V'}^{\geq r}  & \:\: &  \agb^{\leq r} & = & \bigcup_{V'\in\mathfrak{C}^{\leq r}} & Z_{V'}^{\leq r} \\
\agdb^r & = & \bigcup_{V'\in\mathfrak{C}^r_{\mathscr{D}}} & Z_{V',\mathscr{D}}^r & \:\: &  \agdb^{\geq r} & = & \bigcup_{V'\in\mathfrak{C}^{\geq r}_\mathscr{D}} & Z_{V',\mathscr{D}}^{\geq r}  & \:\: & \agdb^{\leq r} & = & \bigcup_{V'\in\mathfrak{C}^{\leq r}_\mathscr{D}} & Z_{V',\mathscr{D}}^{\leq r} \\
\end{array}$$
On obtient ainsi un diagramme commutatif
$$\xymatrix{
\agdb^{\leq r} \ar[d]^{\pib} & \agdb^r \ar[r]\ar[l]\ar[d]^{\pib} & \agdb^{\geq r} \ar[d]^{\pib} \\
\agb^{\leq r} & \agb^r \ar[r]\ar[l] & \agb^{\geq r} }
$$
Dans la proposition suivante, nous explicitons la structure locale de chacun des sch\'emas d\'efinis plus haut. La d\'emonstration est imm\'ediate \`a partir de la proposition~\ref{prop_struc_locale} et de la construction.
 
\begin{proposition}\label{proposition_struc_locale_strates} Pour tout~$V'\in\mathfrak{C}^r$, il existe un isomorphisme canonique stratifi\'e
$$\left(\mvb^r\right) ^{h_{Y_{V'}^r}} / \: \Gamma_{V'} \: \isolong \: \left( \agb^r \right)^{h_{Z_{V'}^r}}\: .$$
Il en est de m\^eme pour $\agb^{\geq r}$, $\agb^{\leq r}$, $\agdb^{r}$, $\agdb^{\geq r}$ et~$\agdb^{\leq r}$. 
\end{proposition} 

Nous avons montr\'e pr\'ec\'edemment que~$\av^r$ \'etait ferm\'e dans~$\av^{\geq r}$. En utilisant la proposition~\ref{proposition_struc_locale_strates}, nous en d\'eduisons le corollaire suivant.

\begin{corollaire} \label{corollaire_ouvert_ferme} Les immersions $\agb^r\hookrightarrow \agb^{\geq r}$, 
$\agdb^r\hookrightarrow \agdb^{\geq r}$, $\agb^{\leq r}\hookrightarrow \agb\times \Spec(\F_p)$ et $\agdb^{\leq r}\hookrightarrow \agdb\times \Spec(\F_p)$ sont ferm\'ees. Les immersions $\agb^r\hookrightarrow \agb^{\leq r}$, 
$\agdb^r\hookrightarrow \agdb^{\leq r}$, $\agb^{\geq r}\hookrightarrow\nolinebreak \agb\times \Spec(\F_p)$ et $\agdb^{\geq r}\hookrightarrow \agdb\times \Spec(\F_p)$ sont ouvertes.
\end{corollaire}

Remarquons que $\mvb^r$ est le sous-sch\'ema de~$\mvb\times\Spec(\F_p)$ o\`u le groupe~$\tilde{G}_0[p]$ est de rang multiplicatif~$p^r$. En effet il existe sur~$\mvb\times\Spec(\F_p)$ une suite exacte $0\rightarrow T_0[p] \rightarrow \tilde{G}_0[p] \rightarrow A_0[p]\rightarrow 0$ dans laquelle le groupe $T_0[p]$ est de rang multiplicatif~$p^{\mathrm{rg}(V')}$, et dans laquelle~$A_0[p]$ est de rang multiplicatif~$p^{r-\mathrm{rg}(V')}$ exactement sur~$\mvb^r$. De m\^eme, $\mvdb^r$ est le sous-sch\'ema de~$\mvdb\times\Spec(\F_p)$ o\`u~$\mathrm{Ker}(\tilde{G}_0\rightarrow\tilde{G}_1)$ est de rang~$p^r$ et \'egal \`a~$\tilde{G}_0[p]^{\mathrm{mul}}$.

Rappelons que d'apr\`es~\cite[th. IV.5.7]{Deg@FaltingsChai}, le sch\'ema ab\'elien universel~$G_0$ sur~$\ag$ s'\'etend en un sch\'ema semi-ab\'elien sur~$\agb$ que l'on note encore~$G_0$. On \'etend le sous-groupe universel~$H$ \`a~$\agdb$ en prenant son adh\'erence dans le groupe quasi-fini et plat~$G_0[p]$ ; on note encore~$H$ le groupe quasi-fini et plat obtenu sur~$\agdb$. Pour tout~$V'\in\mathfrak{C}$, la restriction de~$G_0$ \`a~$Z_{V'}$ correspond \`a la restriction de $\tilde{G}_0$ \`a $Y_{V'}$ par les isomorphismes de la proposition~\ref{prop_struc_locale}. De m\^eme, la restriction de~$H$ \`a~$Z_{V',\mathscr{D}}$ correspond \`a la restriction de~$\mathrm{Ker}(\tilde{G}_0\rightarrow\tilde{G}_1)$ \`a~$Y_{V',\mathscr{D}}$.
La proposition~\ref{proposition_struc_locale_strates} montre par cons\'equent que~$\agb^r$ est le sous-sch\'ema de~$\agb\times\Spec(\F_p)$ o\`u le rang multiplicatif de~$G_0[p]$ est~$p^r$ ; nous avions utilis\'e cette propri\'et\'e pour d\'efinir~$\agb^r$ dans l'introduction. De m\^eme,~$\agdb^r$ est le sous-sch\'ema de~$\agdb\times\Spec(\F_p)$ o\`u~$H$ est multiplicatif de rang~$p^r$ et \'egal \`a~$G_0[p]^{\mathrm{mul}}$. En particulier, $H$ est fini et plat sur~$\agdb^r$.

\subsection{Oubli du niveau} Nous montrons que le morphisme d'oubli du niveau~$\pib$ se comporte exactement de la m\^eme mani\`ere qu'avant compactification.

\begin{proposition} \label{proposition_pib_etale} Le morphisme~$\pib:\agdb\rightarrow\agb$ est \'etale en tout point de~$\agdb^{\geq r}$ et induit une surjection de~$\agdb^{\geq r}$ sur~$\agb^{\geq r}$.
\end{proposition}

\begin{demo} L'\'enonc\'e de surjectivit\'e est clair. On d\'eduit de la proposition~\ref{proposition_struc_locale_strates} que pour montrer que $$\pib \: :\: \agdb\:\longrightarrow\:\agb$$ est \'etale en tout point de~$\agdb^{\geq r}$, il suffit de montrer que pour tout~$V'\in\mathfrak{C}^{\geq r}_{\mathscr{D}}$, le morphisme $$\pib_{V'}\: :\: \mvdb\:\longrightarrow \: \mvb$$ est \'etale en tout point de~$\mvdb^{\geq r}$. Nous allons d\'emontrer ceci par d\'evissage.
D'apr\`es la proposition~\ref{prop_pi_etale_ouvert}, le morphisme naturel
$$\avd\:\longrightarrow\:\av$$
est \'etale en tout point de~$\avd^{\geq r}$. On en d\'eduit que la fl\`eche $\bv\times_{\av} \avd\rightarrow \bv$ est \'etale en tout point de~$\bv^{\geq r}\times_{\av^{\geq r}} \avd^{\geq r}$. Puis le morphisme naturel
$$\bvd \: \longrightarrow \: \bv\times_{\av} \avd$$
est \'etale en tout point de~$\bvd^{\geq r}$. En effet, $\bvd$ param\`etre les diagrammes
\begin{equation}\label{preuve_etale_diagramme1}
\xymatrix{\mathcal{V}_\mathscr{D}^0 / V'^\perp \ar[r]\ar[d]^{\sim} & A_0 \ar[d]^f \\
\mathcal{V}_\mathscr{D}^1 / V'^\perp \ar[r]\ar[d] & A_1 \ar[d] \\
\mathcal{V}_\mathscr{D}^2 / V'^\perp \ar[r]\ar[d] & A_1^t \ar[d]^{f^t} \\
\mathcal{V}_\mathscr{D}^3 / V'^\perp \ar[r] & A_0^t}
\end{equation}
donc le morphisme $\bvd\rightarrow \bv\times_{\av} \avd$ est un~$\mathrm{Ker}(f^t)$-torseur, et $f^t$ est \'etale en tout point de~$\avd^{\geq r}$ par d\'efinition. On en d\'eduit que la fl\`eche $$\mv\times_{\bv} \bvd\:\longrightarrow\: \mv$$ est \'etale en tout point de~$\mv^{\geq r}\times_{\bv^{\geq r}} \bvd^{\geq r}$. Il est \'equivalent de montrer que le morphisme $$\mvd\:\longrightarrow\: \mv$$ est \'etale, de montrer que $\mvd\iso \mv\times_{\bv} \bvd$, et de montrer que les tores relatifs aux torseurs~$\mvd\rightarrow\nolinebreak\bvd$ et~$\mv\rightarrow\bv$ sont les m\^emes, \`a savoir~$\Hom(\mathrm{Sym}^2(V/V'^\perp),\Gm)$. Cela r\'esulte du fait que~$\mvd$ param\`etre les diagrammes
$$\xymatrix{\mathcal{V}_\mathscr{D}^0 / V'^\perp \ar[r]\ar[d]^{\sim} & \tilde{G}_0 \ar[d] \\
\mathcal{V}_\mathscr{D}^1 / V'^\perp \ar[r]\ar[d] & \tilde{G}_1 \ar[d] \\
\mathcal{V}_\mathscr{D}^2 / V'^\perp \ar[r]\ar[d] & \tilde{G}_1' \ar[d] \\
\mathcal{V}_\mathscr{D}^3 / V'^\perp \ar[r] & \tilde{G}_0'}$$
qui rel\`event le diagramme~\ref{preuve_etale_diagramme1}, que sur $\mv\times_{\bv} \bvd$ on conna\^it $\mathcal{V}_\mathscr{D}^0 / V'^\perp \rightarrow \tilde{G}_0$, $\mathcal{V}_\mathscr{D}^3 / V'^\perp \rightarrow \nolinebreak\tilde{G}_0'$ et $\tilde{G}_0\rightarrow\tilde{G}_1\rightarrow\tilde{G}_1'\rightarrow\nolinebreak\tilde{G}_0'$, et que par hypoth\`ese, les $1$-motifs $[\mathcal{V}_\mathscr{D}^1 / V'^\perp \rightarrow \tilde{G}_1]$ et $[\mathcal{V}_\mathscr{D}^2 / V'^\perp \rightarrow \tilde{G}_1']$ doivent \^etre duaux. Il reste \`a prouver que le morphisme $$\mvdb\:\longrightarrow\: \mvb$$ est \'etale en tout point de~$\mvdb^{\geq r}$ ; cela est \'evident car les donn\'ees combinatoires et les structures enti\`eres utilis\'ees pour d\'efinir ces deux plongements toro\"idaux sont les m\^emes, \`a savoir~$\mathscr{S}|_{C(V/V'^\perp)}$ et~$\mathrm{Sym}^2(V/V'^\perp)$.
\end{demo}

\begin{corollaire} \label{corollaire_pib_iso} Le morphisme $\pib:\agdb^r\rightarrow\agb^r$ est un isomorphisme.
\end{corollaire}

\begin{demo} D'apr\`es la proposition~\ref{proposition_pib_etale}, $\pib$ est \'etale et surjectif. Il est en particulier quasi-fini et plat, et le cardinal de ses fibres g\'eom\'etriques d\'ecro\^it par sp\'ecialisation. Comme $\pi:\agd^r\rightarrow\ag^r$ est universellement injectif d'apr\`es le lemme~\ref{lemme_homeo_univ}, on voit que les fibres g\'eom\'etriques~$\pib$ sont des singletons : $\pib$ est un hom\'eomorphisme universel. C'est donc un isomorphisme puisqu'il est \'etale.
\end{demo}

\begin{remarque} Montrer directement l'universelle injectivit\'e de $\pib:\agdb^r\rightarrow\agb^r$ aurait pos\'e un petit probl\`eme de pr\'esentation, car cette question ne se r\'eduit pas imm\'ediatement \`a l'\'etude de $\mvdb^r\rightarrow\mvb^r$ : il faut tenir compte de l'action de~$\Gamma$ et de~$\Gamma_\mathscr{D}$. Nous avons contourn\'e ce probl\`eme en montrant d'abord le caract\`ere \'etale, puis en utilisant la densit\'e de~$\agd^r$ dans~$\agdb^r$.
\end{remarque}

\subsection{Un sous-champ ouvert}
\label{subsection_ouvert}

Dans ce paragraphe, nous construisons un sous-champ ouvert $\agdt$ de $\agdb$ sur $\Spec(\Z)$ tel que le sous-groupe universel $H_0$ soit fini et plat de rang $p^r$ sur $\agdt$ et tel que
$$\agdb^{r} \hookrightarrow \agdt \times \Spec(\F_p)\: .$$
Cette construction nous sera utile pour vérifier que le sous-groupe canonique partiel est fini et plat sur des voisinages stricts convenables.
Notons donc
$$\agdt \: = \: \bigcup_{V'\in\mathfrak{C}^r_\mathscr{D}} \: Z_{V',\mathscr{D}}$$
et rappelons que~$\mathfrak{C}_\mathscr{D}^r$ est l'ensemble des~$V'\in\mathfrak{C}$ tels que~$\mathrm{rg}(V')\leq r$ et $$\mathcal{V}_\mathscr{D}^0/V'^\perp\:=\:\mathcal{V}_\mathscr{D}^1/V'^\perp\:=\:\mathcal{V}_\mathscr{D}^2/V'^\perp\: .$$
D'après les définitions, il est clair que $\agdb^{r} \hookrightarrow \agdt \times \Spec(\F_p)$. Si~$V''$ et~$V'$ sont deux \'el\'ements de~$\mathfrak{C}$, on sait que~$Z_{V',\mathscr{D}}$ est inclus dans l'adh\'erence de~$Z_{V'',\mathscr{D}}$ dans~$\agdb$ si et seulement si~$V''$ est inclus dans~$\gamma\cdot V'$ pour un \'el\'ement~$\gamma$ de~$\Gamma_\mathscr{D}$. Nous en d\'eduisons que~$\agdt$ est un sous-sch\'ema ouvert de~$\agdb$. De plus, pour tout~$V'\in\mathfrak{C}$, la restriction de~$H_0$ \`a~$Z_{V',\mathscr{D}}$ est finie et plate de rang
$$p^{r-\mathrm{rg}(\mathrm{Ker}(\mathcal{V}_\mathscr{D}^0/V'^\perp+pV\:\longrightarrow\:\mathcal{V}_\mathscr{D}^1/V'^\perp+pV))}\: .$$
Ainsi, $H_0$ est fini et plat de rang~$p^r$ sur~$\agdt$.

\section{Construction de sous-groupes canoniques partiels}

\subsection{\'Enonc\'e du th\'eor\`eme principal}
Pour tout sch\'ema~$X$ propre sur~$\Spec(\Z_p)$, on note~$X^{\mathrm{rig}}$ la vari\'et\'e analytique rigide associ\'ee, et pour tout sous-sch\'ema localement ferm\'e~$Y$ de~$X\times\Spec(\F_p)$, on note~$]Y[$ son tube dans~$X^\mathrm{rig}$ ; s'il y a ambiguit\'e sur~$X$, on le notera plut\^ot~$]Y[_X$.

La théorie des champs en géométrie analytique rigide existe sûrement, mais n'a jamais été écrite et n'est vraisemblablement pas complètement triviale. Comme nous l'avons mentionné dans l'introduction, nous rajoutons donc des structures de niveau principales en un entier $n\geq 3$ premier à $p$ à tous les champs $\ag$, $\agb$, $\agb^r$, $\cdots$, $\agd$, $\agdb$, $\agdb^r$, $\cdots$. Cela revient juste à considérer les produits fibrés des objets de la liste précédente par le morphisme d'oubli $\agn \rightarrow \ag$ sur $\Spec(\Z[1/n])$. Les morphismes étales restent étales, et les isomorphes des isomorphismes. Mais les champs algébriques sont transformés en schémas, ce qui nous permet d'étudier les variétés analytiques rigides associées.

D'apr\`es le corollaire~\ref{corollaire_pib_iso} et la remarque suivant le corollaire~\ref{corollaire_ouvert_ferme}, la formule $H^r_\mathrm{mul} =\nolinebreak (\pib^{-1})^* H_0$
d\'efinit un groupe multiplicatif fini, plat et totalement isotrope de rang~$p^r$ sur~$\agnb^r$. Le groupe~$H^r_\mathrm{mul}$ est inclus dans le groupe quasi-fini et plat~$G_0[p]$ et v\'erifie $$H^r_{\mathrm{mul},s}=G_{0,s}[p]^\mathrm{mul}$$ pour tout~$s\in\agnb^r$. Le th\'eor\`eme de rigidit\'e des groupes de type multiplicatif~\sga{3}{iv}{th. 3.bis} montre que~$H^r_\mathrm{mul}$ se rel\`eve canoniquement en un sous-groupe~$H^r_\mathrm{can}\subset G[p]$ qui est fini, plat et totalement isotrope de rang~$p^r$ sur~$]\agnb^r[$.
Nous nous demandons s'il existe un voisinage strict sur lequel~$H^r_\mathrm{can}$ se prolonge. Le th\'eor\`eme suivant r\'epond affirmativement \`a cette question.

\begin{theoreme} \label{theoreme_principal} Il existe un voisinage strict~${U}$ de~$]\agnb^r[$ dans $]\agnb^{\leq r}[$ et un sous-groupe~$H^r_\mathrm{can,{U}}$ de~$G_0[p]$ fini et plat de rang~$p^r$ sur~${U}$ dont la restriction \`a $$]\agnb^r[$$ est \'egale \`a~$H^r_\mathrm{can}$.
Le sous-groupe~$H^r_\mathrm{can,{U}}$ est unique et totalement isotrope d\`es que toute composante connexe de~$U$ rencontre $$]\agnb^{r}[\: .$$
\end{theoreme}

Le groupe~$H^r_\mathrm{can,{U}}$ est appell\'e \og sous-groupe canonique partiel de rang~$p^r$ sur~${U}$ \fg. Sa restriction \`a l'union des composantes connexes de~$U$ qui rencontrent $$]\agnb^r[$$ est canonique. Remarquons qu'il r\'esulte de~\cite[d\'ef. 1.2.1]{Rigide@Berthelot} que cette union disjointe forme encore un voisinage strict de~$]\agnb^r[$.

\begin{remarque} Dans le cas o\`u~$r=g$, le sch\'ema~$\agnb^{g}$ est le lieu ordinaire de~$\agnb^{\leq g}=\nolinebreak\agnb$. Le groupe~$H^g_\mathrm{can}$ rel\`eve \`a $$]\agnb^g[$$ le noyau du morphisme de Frobenius de~$G_0$ sur~$\agnb^g$, et les sous-groupes canoniques partiels de rang~$p^g$ sont les sous-groupes canoniques habituels. La d\'emonstration du th\'eor\`eme~\ref{theoreme_principal} fournit dans ce cas une nouvelle preuve --- non-effective --- de certains r\'esultats de~\cite{Canonique@AbbesMokrane} et de~\cite{Canonique@Andreatta}.
\end{remarque}

La preuve du th\'eor\`eme~\ref{theoreme_principal} occupe les deux paragraphes suivants. Nous y d\'emontrons respectivement l'unicit\'e puis l'existence.

\subsection{D\'emonstration de l'unicit\'e} Soit~$U$ un voisinage strict de~$]\agnb^r[$ dont chaque composante connexe rencontre $$]\agnb^r[\: .$$
Montrons qu'il existe au plus un sous-groupe canonique partiel de rang~$p^r$ sur~$U$. Notons
$\agn^{\mathrm{an}}$ et $\agnd^{\mathrm{an}}$
les champs alg\'ebriques rigides localement de type fini obtenus par analytification de $$\agn\times\Spec(\Q_p)$$ et de $$\agnd\times\Spec(\Q_p)\: .$$
L'ouvert $\agn^{\mathrm{an}}$ de~$\agnb^{\mathrm{rig}}$ est donc Zariski-dense. Rappelons le résultat suivant de géométrie rigide~\cite[th. 1.6 I)]{Rigid@Lu}.

\begin{proposition}\label{prop_Bartenwerfer} Soit~$X$ une variété rigide lisse et~$Z$ un fermé de Zariski de~$X$. Notons~$V$ l'ouvert complémentaire de~$Z$ dans~$X$. La variété rigide~$V$ est irréductible.
\end{proposition}

\begin{corollaire} \label{coro_stupide} Toutes les composantes connexes de~$U\cap \agn^{\mathrm{an}}$ rencontrent~$\agnb^r\cap\agn^{\mathrm{an}}$.
\end{corollaire}

\begin{demo} L'ouvert~$U$ est réunion de ses composantes connexes, qui sont irréductibles puisque $$\agnb^{\mathrm{rig}}$$ l'est. Les composantes connexes de~$U\cap\agn^{\mathrm{an}}$ sont en bijection avec celles de~$U$ d'après la proposition~\ref{prop_Bartenwerfer}. Comme toutes les composantes connexes de~$U$ rencontrent $$]\agnb^r[\: ,$$
il en est de même des composantes connexes de~$U\cap \agn^{\mathrm{an}}$.
\end{demo}

Nous pouvons \`a pr\'esent d\'emontrer qu'il existe au plus un sous-groupe canonique partiel de rang~$p^r$ sur~$U$. Soient~$H$ et~$H'$ deux tels sous-groupes. Leur restriction \`a~$U\cap\agn^\mathrm{rig}$ d\'eterminent deux sections~$s$ et~$s'$ de $$\pi:\agnd^{\mathrm{rig}}\rightarrow\agn^{\mathrm{rig}}$$ sur~$U\cap\agn^\mathrm{rig}$ qui co\"incident sur~$]\agnb^r[\cap\agn^\mathrm{rig}$. D'apr\`es le principe de prolongement analytique~\cite[prop.~0.1.13]{Rigide@Berthelot} et le corollaire~\ref{coro_stupide}, les sections~$s$ et~$s'$ co\"incident sur la vari\'et\'e int\`egre~$U\cap\agn^\mathrm{rig}$. Ainsi, les restrictions de~$H$ et~$H'$ \`a $$U\cap\agn^\mathrm{rig}$$ sont \'egales. Comme~$U\cap\agn^\mathrm{rig}$ est dense dans~$U$ et~$H$ (resp. $H'$) est fini et plat sur~$U$, le groupe~$H$ (resp.~$H'$) est \'egal \`a l'adh\'erence de~$H_U$ (resp. $H'_U$) dans~$G_0$. On a donc~$H=H'$ sur~$U$, ce qui termine la d\'emonstration de l'unicit\'e.

\subsection{D\'emonstration de l'existence}
Dans~\cite[th.~1.3.5]{Rigide@Berthelot}, Berthelot \'etudie le probl\`eme g\'en\'eral suivant de g\'eom\'etrie rigide. Soit
$$\xymatrix{
 &  Y' \ar[r]^{i'} \ar[d]^v & P' \ar[d]^u \\
X \ar[ru]^{j'}  \ar[r]^j & Y \ar[r]^i & P}$$
un diagramme commutatif o\`u~$P$ et~$P'$ sont des sch\'emas formels de type fini sur~$\mathrm{Spf}(\Z_p)$, o\`u $Y$,~$Y'$ et~$X$ sont des sch\'emas de type fini sur~$\Spec(\F_p)$, o\`u~$j$ et~$j'$ sont des immersions ouvertes, et o\`u~$i$ et~$i'$ sont des immersions ferm\'ees. On suppose que~$u$ est propre et que~$v$ est \'etale sur un voisinage ouvert de~$X$ dans~$P'$. Berthelot montre alors que~$u$ induit un isomorphisme entre un voisinage strict de~$]X[_{P'}$ dans~$]Y'[_{P'}$ et un voisinage strict de~$]X[_{P}$ dans~$]Y[_{P}$.

Notons~$\agnb^\mathrm{for}$ et~$\agndb^\mathrm{for}$ les compl\'etions formelles~$p$-adiques de~$\agnb$ et de~$\agndb$, puis appliquons le r\'esultat de Berthelot au diagramme commutatif
$$\xymatrix{
\agndb^r\ar[d]^{\pib}_{\sim} \ar[r] &  \agndb^{\leq r} \ar[r] \ar[d]^{\pib} & \agndb^\mathrm{for} \ar[d]^{\pib} \\
\agnb^r \ar[r] & \agnb^{\leq r} \ar[r] & \agnb^\mathrm{for}}$$
dans lequel~$\pib$ est propre, r\'ealise un isomorphisme entre~$\agndb^r$ et~$\agnb^r$ d'apr\`es le corollaire~\ref{corollaire_pib_iso}, et est \'etale sur un voisinage ouvert de $$\agndb^{\geq r}$$ dans $\agndb^\mathrm{for}$  d'apr\`es la proposition~\ref{proposition_pib_etale} et le caract\`ere quasi-compacit\'e de~$\agndb$. Nous en d\'eduisons l'existence d'un voisinage strict~$V$ de $$]\agndb^r[$$ dans~$]\agndb^{\leq r}[$ et d'un voisinage strict~$U$ de~$]\agnb^r[$ dans~$]\agnb^{\leq r}[$ tel que~$\pib$ r\'ealise un isomorphisme entre~$V$ et~$U$. Le groupe~$(\pib^{-1})^* H_0$ est alors quasi-fini et plat sur~$U$, et sa restriction \`a $$]\agnb^r[$$ est \'egale \`a~$H^r_\mathrm{can}$. Nous aurons d\'emontr\'e le th\'eor\`eme~\ref{theoreme_principal} si, quitte \`a rapetisser~$U$, nous prouvons que $(\pib^{-1})^* H_0$ est fini et plat de rang~$p^r$ sur~$U$.

\begin{remarque} Il est clair que~$(\pib^{-1})^* H_0$ est fini et plat de rang~$p^r$ sur le lieu~$U\cap\: ]\agn[$ des vari\'et\'es ab\'eliennes ayant bonne r\'eduction.
\end{remarque}

Pour montrer que, quitte \`a rapetisser~$U$, le groupe $(\pib^{-1})^* H_0$ est fini et plat de rang~$p^r$ sur~$U$, il suffit de montrer que, quitte \`a rapetisser~$V$, le groupe~$H_0$ est fini et plat de rang~$p^r$ sur~$V$. Notons $\agndt$ le schéma sur $\Spec(\Z[1/n])$ obtenu par changement de base de $\agdt$ \textit{via} $\agn\rightarrow \ag$ (\textit{cf}. paragraphe~\ref{subsection_ouvert}).
 Notons $$\agndt^{\mathrm{\: rig}}$$ la vari\'et\'e analytique localement de type fini associ\'ee au sch\'ema~$\agndt\times \Spec(\Q_p)$. Elle forme un ouvert de Zariski dans $$\agndb^{\mathrm{rig}}$$ et~$H_0$ est fini et plat de rang~$p^r$ sur~$\agndt^{\mathrm{\: rig}}$. Comme~$\agndb^r$ est inclus dans~$\agndt\times \Spec(\F_p)$, tout voisinage strict de $$]\agndb^r[$$ dans $]\agndb^{\leq r}[$ contient un voisinage strict inclus dans~$]\agndb^{\leq r}[\:\cap\:\agndt^{\mathrm{\: rig}}$. Quitte \`a rapetisser~$V$, on peut donc supposer que~$H_0$ est fini et plat de rang~$p^r$ sur~$V$, ce qui termine la d\'emonstration du th\'eor\`eme~\ref{theoreme_principal}.

\textbf{Remerciements.} {Nous remercions le rapporteur pour sa relecture attentive de cet article.}

\providecommand{\bysame}{\leavevmode ---\ }
\providecommand{\og}{}
\providecommand{\fg}{}
\providecommand{\smfandname}{et}
\providecommand{\smfedsname}{\'eds.}
\providecommand{\smfedname}{\'ed.}
\providecommand{\smfmastersthesisname}{M\'emoire}
\providecommand{\smfphdthesisname}{Th\`ese}

\end{document}